\documentclass[a4paper,12pt]{article}
\usepackage{amssymb,amsmath}

\newcommand{\F}{\mathbb{F}}
\renewcommand{\b}[1]{{\bf #1}}
\newcommand{\cl}{\ell}
\newcommand{\cN}{\mathcal{N}}

\newcommand{\Bf}{\mathbf{f}}

\newcommand{\x}{\mathbf{x}}

\newcommand{\Af}{\mathbb{A}^n(\F_q)}
\newcommand{\Aft}{\mathbb{A}^t(\F_q)}

\renewcommand{\mod}[1]{\hspace{-2.9mm}\pmod{#1}}

\newtheorem{theorem}{Theorem}
\newtheorem{coroll}{Corollary}

\newtheorem{lemma}{Lemma}

\begin{document}
\title{A Note on the Chevalley--Warning Theorems}
\author{D.R. Heath-Brown\\Mathematical Institute, Oxford}
\date{}
\maketitle
\begin{center}
{\em In memory of Anatolii Karatsuba}
\end{center}

\section{Introduction}

The most widely known result of Chevalley--Warning type states that if
one has a polynomial over a finite field of characteristic $p$, and
the number of variables exceeds the degree, then the number of zeros
is a multiple of $p$.  In particular, if the polynomial is
homogeneous, there is at least one non-trivial zero.

There are a number of related results in the literature.
Generally, let $\F_q$ be a finite field, and let $\Bf=(f_1(\x),\ldots,f_r(\x))$
be an $r$-tuple of polynomials $f_i(x_1,\ldots,x_n)\in\F_q[x_1,\ldots,x_n]$.
Let $d_i$ be the total degree of $f_i$, and write $d=d_1+\ldots+d_r$.
For any subset $S\subseteq\F_q^n$ we put
\[z(\Bf;S)=z(S):=\{\x\in S:\Bf(\x)=0\}\]
and
\[\cN(\Bf;S)=\cN(S):=\# z(\Bf;S).\]
Then the original 1935 result of Chevalley \cite{Chev} stated that if
$n>d$ and $z(\F_q^n)$ is non-empty, then it contains at least 2 points.
Immediately
afterwards it was shown by Warning \cite{Warn} that if $n>d$ then 
\begin{equation}\label{chev}
p\mid \cN(\F_q^n), 
\end{equation}
where $p$ is the characteristic of $\F_q$.  In fact this followed
from a trivial re-arrangement of Chevalley's argument. In the same
paper Warning proved that if $n\ge d$, and if $H_1$
and $H_2$ are parallel affine hyperplanes in $\F_q^n$ then
\begin{equation}\label{war}
\cN(H_1)\equiv \cN(H_2)\mod{p}.  
\end{equation}
As a corollary he deduced that if $z(\F_q^n)$ is non-empty then
\begin{equation}\label{war2}
\cN(\F_q^n)\ge q^{n-d}\;\;\;(n>d).
\end{equation}
Later, in 1964, Ax \cite{Ax}
strengthened (\ref{chev}) by showing
that indeed one has 
\begin{equation}\label{ax}
q\mid \cN(\F_q^n)\;\;\;(n>d).
\end{equation}
When $n$ is large compared with $d$ there are further improvements
possible, as was shown by Katz \cite{K} for example. These sharpenings
of (\ref{ax}) lead to improvements of all the results we will
establish.  However we wish to focus
here on the situation in which $n$ is not much larger than $d$, and
so will be content with using (\ref{ax}).

Our first result shows that there is a strengthening of (\ref{war})
corresponding to Ax's improvement of (\ref{chev}). To stress that the
affine linear spaces we encounter are not necessarily vector subspaces
of $\F_q^n$ we will work over $\Af$, but of course the sets $z(\Af)$
and $z(\F_q^n)$ are the same.
\begin{theorem}\label{affax}
With the notation above, we have 
\begin{equation}\label{affaxeq}
\cN(L_1)\equiv \cN(L_2)\mod{q}
\end{equation}
for any two parallel linear spaces $L_1,L_2\subseteq\Af$ of
dimension $d$ or more.
\end{theorem}
Indeed our argument shows how one may deduce (\ref{affaxeq}) from
(\ref{ax}).  

We next look at improvements to (\ref{war2}).  If
$K=\mathbb{F}_{q^k}$ has field basis $\omega_1,\ldots,\omega_k$ over
$F=\F_q$ we set
\begin{equation}\label{norm}
N_k(x_1,\ldots,x_k):=N_{K/F}(x_1\omega_1+\ldots+x_k\omega_k).
\end{equation}
This produces a form of degree $k$ in $k$ variables, with the property
that $N_k(x_1,\ldots,x_k)=0$ with variables $x_i\in\ F_q$ only when
the $x_i$ are all zero.  In particular, if $\Bf$ consists of a single 
polynomial $f_1=N_d$ then 
$z(\Bf;\Af)$ will be a linear space of dimension $n-d$, and we will
have equality in (\ref{war2}).  In contrast we have the following
results.
\begin{theorem}\label{lin}
Suppose that $n>d$ and that $z(\Af)$ is non-empty, and is not a
linear subspace of $\Af$.  Then
\begin{enumerate}
\item[(i)] For any $q$ we have $\cN(\Af)>q^{n-d}$;
\item[(ii)] If $q\ge 4$ we have $\cN(\Af)\ge 2q^{n-d}$; and
\item[(iii)] For any $q$ we have $\cN(\Af)\ge q^{n+1-d}/(n+2-d)$
  providing that the polynomials $\Bf$ are homogeneous.
\end{enumerate}
\end{theorem}
Our proof of this will use only combinatorial facts about $\Af$, along with
(\ref{affaxeq}) and its consequence (\ref{ax}).  Thus the analogues 
of these results would be
true if one replaced $z(\Af)$ by any set with the property 
(\ref{affaxeq}) for every pair of 
parallel $d$-dimensional linear spaces $\Lambda_1,\Lambda_2\subseteq\Af$.
In particular we see that part (ii) is best possible in this sense,
since (\ref{affaxeq}) would hold if $z(\Af)$ were composed of two
parallel linear spaces of dimension $n-d$.

It may be instructive to examine two simple examples at this point.
The first shows that one cannot hope to
improve part (iii) of Theorem \ref{lin} to say that $\cN(\Af)\ge q^{n+1-d}$.
Here we take $\Bf$ to consist of the single polynomial
$Q(x_1,\ldots,x_4)$, where 
\[Q(x_1,x_2,x_3,x_4)=x_1x_2+x_3^2+x_3x_4+cx_4^2\]
with the coefficient $c\in\F_q$ chosen so that
$x_3^3+x_3x_4+cx_4^2$ does not factor over $\F_q$.
Then $Q=0$ has precisely $q^3-q^2+q$ solutions
$(x_1,\ldots,x_4)\in\F_q^4$, and since $d=2$ in this case we have 
\[\cN(\Bf;\Af)=q^{n-4}(q^3-q^2+q)=q^{n+1-d}(1-q^{-1}+q^{-2})<q^{n+1-d}. \]
Indeed this also 
shows that part (ii) of the theorem cannot be extended to $q=2$.  

Our second example shows that when $\Bf$ consists of a single polynomial
$f_1$, it is possible for $z(\Af)$ to be a linear
space of dimension $n-k$ even though $f_1$ does not split into linear
factors. For this we take $n=4$ and choose a field basis
$1,\alpha$ for $\F_{q^2}$ over $\F_q$.  Let $N$ and $N'$ denote norm
forms for $\F_{q^4}$ over $\F_{q^2}$  and $\F_{q^2}$ over $\F_{q}$
respectively, using the basis $1,\alpha$ in the latter case.  Let
$\sigma$ be the nontrivial automorphism of $\F_{q^2}$ over
$\F_q$. Then we may write
$N(\x)=Q_1(\x)+\alpha Q_2(\x)$ for certain quadratic forms
$Q_1,Q_2$ over $\F_q$.  It is clear that we may use a linear change of
variables over $\overline{\F_q}$ to write $Q_1+\alpha Q_2=Y_1Y_2$ and
$Q_1+\alpha^{\sigma} Q_2=Y_3Y_4$ with independent variables
$Y_1,\ldots,Y_4$.  It follows that if we take
$\beta\in\F_{q^2}\setminus\F_q$ different from $\alpha$ and
$\alpha^{\sigma}$ then $Q_1+\beta Q_2=$ will have rank 4.  Such
a $\beta$ certainly exists, except when $q=2$.  We now see that
the polynomial $f(\x)=(Q_1(\x)+\beta Q_2(\x))
(Q_1(\x)+\beta Q_2(\x))^{\sigma}$ is defined over $\F_q$ and does not 
split into linear factors.  None the less we have $f=0$ only for $\x=0$.

Finally we will investigate the structure of $z(\Bf;\Af)$ geometrically.
\begin{theorem}\label{geom}
Given $n$ and $d$ there is a constant $C(n,d)$ as follows.  For any
set of polynomials $\Bf$ whose degrees total $d$ there is a
corresponding set of absolutely irreducible varieties $V_1,\ldots,V_m$
defined over $\F_q$ such that $z(\Bf;\Af)$ is precisely the set of
$\F_q$-points in $V_1\cup\ldots\cup V_m$.  Moreover each variety $V_i$
has degree at most $C(n,d)$, and $m\le C(n,d)$.
\end{theorem}
The Lang-Weil bound \cite{lw} now gives us the following immediate
corollary.
\begin{coroll}\label{cor}
Suppose that the largest dimension among the varieties $V_i$ is $D$,
and that this occurs $k$ times.  Then, enlarging $C(n,d)$ if
necessary, we have
\begin{equation}\label{lw}
|\cN(\Bf;\Af)-kq^D|\le C(n,d)q^{D-1/2}.
\end{equation}
Thus if $q\gg_{n,d}1$ we have $D\ge n-d$, and indeed if $z(\Bf;\Af)$ is not a 
linear space of dimension $n-d$ then $D\ge n+1-d$. In the latter case
we have
\[\cN(\Bf;\Af)\ge q^{n+1-d}-C(n,d)q^{n+1/2-d}\]
for $q\gg_{n,d}1$.
\end{coroll}
One would conjecture that $D\ge n-d$ irrespective of the size of $q$,
and it would be of interest to have a direct geometrical proof of
this.  Of course, the lower bounds on $D$ follow on comparing
(\ref{lw}) with (\ref{war2}) and part (iii) of Theorem \ref{lin}.

{\bf Acknowledgement} This research was begun while the author was
visiting the Institute for Advanced Study, in Princeton.  The
hospitality and financial support of the institute is gratefully
acknowledged. 

Thanks are also due to the anonymous referee for his careful reading
of the original version of this paper, which resulted in the
correction of two errors.

Finally, we would like to record our thanks to Professor David Leep,
who pointed out a difficulty with our original treatment of Theorem 1.
\bigskip

\section{Proofs}

We now present our proof of Theorem \ref{affax}. In view of (4) the
result is trivial if $L_1$ and $L_2$ have dimension strictly greater
than $d$, since we would have $q\mid\cN(L_1)$ and $q\mid\cN(L_2)$ in
this case.  We therefore suppose that $L_1$ and $L_2$ have dimension equal
to $d$.  Let $L=<\b{e}_1,\ldots,\b{e}_d>$ be the linear space parallel
to $L_1$ and $L_2$ but passing through the origin, and let
$L_j=L+\b{c}_j$ for $j=1$ and 2. Then
$\cN(\b{f};L_j)=\cN(\b{g}^{(j)};\mathbb{A}^d(\F_q))$, where
\[g_i^{(j)}(y_1,\ldots,y_d)=f_i(\sum_{k=1}^dy_k\b{e}_k+\b{c}_j)
\;\;\;(j=1,2,\;1\le i\le r).\]
It could happen that the terms of degree $d_i$ in $g_i^{(1)}$ all
vanish, but this happens if and only if the corresponding terms in
$g_i^{(2)}$ also vanish. In this situation the total degree of each of
the systems $(g_1^{(1)},\ldots,g_r^{(1)})$ and
$(g_1^{(2)},\ldots,g_r^{(2)})$ will be strictly less than $d$, so that 
$q\mid\cN(L_1)$ and $q\mid\cN(L_2)$, by (4).  It follows that we may
assume that the leading homogeneous parts of $g_i^{(1)}$  and
$g_i^{(2)}$, both have degree $d_i$.  Indeed we then see that these
leading homogeneous parts are the same.

Given a polynomial $f(x_1,\ldots,x_n)$ of total degree $e$ there are two
reasonable ways of associating a form to it.  One may take
$f_-(x_1,\ldots,x_n)$ to be the homogeneous part of degree $e$, or one
may define 
\[f_+(x_0,\ldots,x_n)=x_0^ef(x_1/x_0,\ldots,x_n/x_0).\]
For a system $\Bf$ we define $\Bf_-$ and $\Bf_+$ by the above
processes, using degree $e=d_i$ for each polynomial $f_i$.  Clearly each
zero of $\Bf$ produces exactly $q-1$ zeros of $\Bf_+$ with
$x_0\not=0$; and the zeros of $\Bf_+$ with $x_0=0$ correspond
precisely to the zeros of $\Bf_-$.  Thus
\[\cN(\Bf_+;\F_q^{n+1})=(q-1)\cN(\Bf;\F_q^n)+\cN(\Bf_-;\F_q^n).\]
In particular, if $n\ge d$ then (\ref{ax}) yields $q\mid 
\cN(\Bf_+;\F_q^{n+1})$ and hence 
\[\cN(\Bf;\F_q^n)\equiv\cN(\Bf_-;\F_q^n)\mod{q}.  \]
We therefore see that the value of $\cN(\Bf;\F_q^n)$ modulo $q$ only 
depends on the leading homogeneous
parts of the polynomials $f_1,\ldots,f_r$.  

Now suppose we have parallel
$d$-dimensional affine linear spaces $L_1$ and $L_2$ in $\F_q^n$.  We
represent the restriction of the polynomials $\Bf$ to $L_1$ by a set of
polynomials $(g^{(1)}_1,\ldots,g^{(1)}_r)$ in $d$ variables, and
similarly for $L_2$.  As discussed above, the leading
homogeneous parts of $g^{(1)}_i$ and $g^{(2)}_i$ will be the same, and
we are therefore led to Theorem \ref{affax}.
\bigskip

Our proof of Theorem \ref{lin} is based on the following two lemmas.
\begin{lemma}\label{l1}
Let $L_0\subseteq\Af$ be a linear space.  Choose
a linear space $L$ of maximal dimension, $k$ say, such that $L\supseteq
L_0$ and $\cN(L)=\cN(L_0)$.  Suppose $L'\supset L$ is a linear space of 
dimension $k+1$ such that $\cN(L')$ is minimal.  Then
\begin{equation}\label{basic}
\cN(\Af)\ge \cN(L)+\frac{q^{n-k}-1}{q-1}(\cN(L')-\cN(L)).
\end{equation}
\end{lemma}
\begin{lemma}\label{l2}
Let $S\subseteq\Aft$ be a set containing $t+1$ points in general position.
Then
\begin{itemize}
\item[(i)]  If $q=2$, and if there is no 2-plane $L\subseteq\Aft$
  meeting $S$ in exactly 3 points, then $S=\Aft$.
\item[(ii)] If $q\ge 3$, and if $\cl\subseteq S$ for every line $\cl$
  meeting $S$ in at least two points, then $S=\Aft$.
\item[(iii)] If $q\ge 4$, and if $\#(S\cap\cl)\ge q-1$ 
for every line $\cl$ meeting $S$ in at least two points, then 
$\Aft\setminus S$ is contained in a hyperplane.
\item[(iv)] If $m\ge 2$ is an integer, and if $\#(S\cap\cl)\ge m+1$ 
for every line $\cl$ meeting $S$ in at least two points, then 
\[\# S\ge\frac{m^{t+1}-1}{m-1}.\]
\end{itemize}
\end{lemma}

To deduce Theorem \ref{lin} we assume that $z(\Af)$ contains a maximal
set of $t\ge n+1-d$
points in general position.  We take $\Aft$ to be
the space spanned by these points and apply Lemma \ref{l2} to
the set $S$ consisting of zeros of $\Bf$ lying in $\Aft$.  We will
give the details required for the different parts of Theorem \ref{lin}
in due course, but the general strategy is as follows. Lemma~\ref{l2} will
either provide a suitable lower bound for $\cN(\Af)$, or produce a
linear space $L_0$ to which we will apply Lemma \ref{l1}.  This linear space
will have dimension two in the case of part (i) of Lemma \ref{l2}, or
dimension one otherwise.  We will have arranged that $\cN(L_0)\le
q-1$, and since $\cN(L)=\cN(L_0)$ in Lemma~\ref{l1} the dimension $k$
of $L$ can be at most $d$, for otherwise (\ref{affaxeq}) will give us a
contradiction.  If $k\le d-2$ then (\ref{basic}) implies that
\[\cN(\Af)\ge \frac{q^{n-k}-1}{q-1}\ge \frac{q^{n+2-d}-1}{q-1}\geq^{n+1-d},\]
which is a satisfactory lower bound.
Thus we may assume that either $L$ or $L'$ has dimension $d$.  Let
$L_1=L$ or $L'$ as appropriate.  If $\cN(L_1)$ takes a value 
$v\le q-1$ then we may
deduce from (\ref{affaxeq}) that $\cN(L^*)\ge v$ for every affine $d$-plane
$L^*$ parallel to $L_1$.  Covering $\Af$ with
such linear spaces we deduce that $\cN(\Af)\ge vq^{n-d}$.  Such bounds
as these will suffice in all cases for the theorem.

We begin by establishing part (i) of Theorem \ref{lin}.  From parts
(i) and (ii) of Lemma \ref{l2} we may obtain a satisfactory bound 
$\cN(\Af)\ge q^{n+1-d}$
unless either $q=2$ and there is a 2-plane $L_0\subseteq \Af$ with
$\cN(L_0)=3$, or $q\ge 3$ and there is a line $L_0\subseteq \Af$ with
$2\le\cN(L_0)\le q-1$.  We now apply Lemma \ref{l1}.  As above we
must have $k\le d$.  Moreover if $k\le d-1$ then (\ref{basic}) yields
\[\cN(\Af)\ge \frac{q^{n+1-d}-1}{q-1}>q^{n-d},\]
which is satisfactory. If $k=d$ and $q=2$ then (\ref{basic})
produces
\[\cN(\Af)\ge 3+\frac{2^{n-d}-1}{2-1}>2^{n-d}.\]
Finally, if $k=d$ and $q\ge 3$ then $2\le\cN(L_0)=\cN(L)\le q-1$, and
the argument above, with $2\le v\le q-1$, shows that 
$\cN(\Af)\ge 2q^{n-d}$.  Thus $\cN(\Af)>q^{n-d}$ in all cases.

We turn now to part (ii) of Theorem \ref{lin}.  By part (iii) of Lemma
\ref{l2} we have $\cN(\Af)\ge q^{n+1-d}-q^{n-d}\ge 2q^{n-d}$ unless
there is a line $L_0\subseteq \Af$ with $2\le\cN(L_0)\le q-2$.  As above,
when we apply Lemma \ref{l1} we may assume that $k=d-1$ or $d$.
If $k=d$ then $2\le \cN(L)\le q-2$ and we may
apply (\ref{affaxeq}) as before to conclude that $\cN(\Af)\ge 2q^{n-d}$.
Similarly if $k=d-1$ and $\cN(L')-\cN(L)=1$ we have
$3\le \cN(L')\le q-1$ and therefore $\cN(\Af)\ge 3q^{n-d}$.  Finally,
if $k=d-1$ and $\cN(L')-\cN(L)\ge 2$ we deduce from (\ref{basic}) that
\[\cN(\Af)\ge 2\frac{q^{n+1-d}-1}{q-1}\ge 2q^{n-d}.\]
This establishes the required bound in all cases.

To prove the third claim of Theorem \ref{lin} we apply Lemma \ref{l2} 
part (iv) with an integer $m\le q-1$ to be chosen in due course.  Thus
\[\cN(\Af)\ge\frac{m^{n+2-d}-1}{m-1}\ge m^{n+1-d}\]
unless there is a line with $2\le\cN(L_0)\le m$.  In the latter case
we deduce from Lemma \ref{l1} that
\[\cN(\Af)\ge \frac{q^{n-k}-1}{q-1}(\cN(L')-m),\]
and by the previous argument it suffices to consider the values
$k=d-1$ and $k=d$. If $k=d-1$ we deduce that $\cN(\Af)\ge
q^{n-d}(q-m)$, unless 
\[m+1\le\cN(L')\le q-1.  \]
However, in the latter
case we conclude by our standard argument using (\ref{affaxeq}) that
$\cN(\Af)\ge (m+1)q^{n-d}$.  In the homogeneous case the value $k=d$
cannot occur.  To show this we consider two cases.  If $k=d$ and 
$\mathbf{0}\in L$ then $q-1\mid \cN(L)-1$, since if $\x\in z(L)$ then
every scalar multiple of $\x$ is also in $z(L)$.  This
however is impossible since $2\le\cN(L_0)=\cN(L)\le m$.  On the other
hand, if $\mathbf{0}\not\in L$ we consider the
$(d+1)$-dimensional linear space $L':=\langle L,\mathbf{0}\rangle$.  Here we
find that $\cN(L')=1+(q-1)\cN(L)$.  According to (\ref{ax}) we will
have $q\mid \cN(L')$, whence $\cN(L)\equiv 1\pmod{q}$.  This again is
impossible, since $2\le\cN(L)\le m$.  It follows that one of the
inequalities \[\cN(\Af)\ge m^{n+1-d}, \]
or 
\[\cN(\Af)\ge q^{n-d}(q-m),\] 
or 
\[\cN(\Af)\ge (m+1)q^{n-d}\]
must hold.  The required estimate now
follows on choosing
\[m=\left[\frac{q(n+1-d)}{n+2-d}\right],\]
and noting that 
\[\left[\frac{qv}{v+1}\right]^v\ge\frac{q^v}{v+1}\]
for positive integers $v\ge 2$ and $q\ge 2v+3$ with the exception of
$v=2$, $q=7$.  Since $q-1\mid\cN(\Af)-1$, 
the remaining cases follow from parts (i) and (ii) of
the theorem, together with (4).
\bigskip

We now prove Lemmas \ref{l1} and \ref{l2}.  For the first of these it
is enough to observe that $\Af$ is the disjoint union of $L$ together
with the sets $L^*\setminus L$, with $L^*$ running over all
$(k+1)$-dimensional linear spaces containing $L$.  There are
$(q^{n-k}-1)/(q-1)$ such spaces $L^*$.  Thus the lemma follows from
the fact that $\cN(L^*\setminus L)\ge\cN(L')-\cN(L)$.

Lemma \ref{l2} will require distinctly more work.  We prove the four
claims separately, using induction on $t$ in each case.  The various
statements are all trivially true for $t=0$ and $t=1$.  To handle the
induction step for case (i) we may assume that $S$ contains a linear
space $L_0$ say, of dimension $t-1$, together with a point $P_0$ say,
not on $L_0$.  Now choose any point $P\not\in L_0$ with $P\not=P_0$,
and aim to show that $P\in S$.  We will then be able to conclude that
$S=\Aft$ as required.  Let $P_1$ be any point in $L_0$ and consider
the 2-plane generated by $P,P_0$ and $P_1$.  Since $q=2$, this plane 
consists of the three generators together with a fourth point $P_2$ 
say, which must belong to $L_0$.  Our construction has ensured that
$P_0,P_1$ and $P_2$ are all in $S$, and so by our hypothesis we also
have $P\in S$, as required.

For part (ii) we use the same initial setup.  We assume that $S$ 
contains a linear space $L_0$ of
dimension $t-1$, along with a point $P_0$ not in $L_0$, and we choose
a point $P\not\in L_0$ different from $P_0$. Suppose firstly that the line $\cl$
generated by $P_0$ and $P$ meets $L_0$, at a point $P_1$ say.  Then
$\ell$ meets $S$ in at least 2 points, namely $P_0$ and $P_1$.  Then,
by our hypothesis, $\cl$ is contained in $S$, whence $P$ in particular
belongs to $S$.  This deals with all points $P$ except those which lie
on the hyperplane, $L_1$ say, which is parallel to $L_0$ and which 
passes through $P_0$.  To
handle such points $P$ we begin by fixing any point $P_1$ on $L_0$.
We then consider the line $\cl$ generated by $P$ and $P_1$.
Since $q\ge 3$ this line contains at least one point additional point
$P_2$ say, which cannot lie in $L_1$.  Thus $P_1\in S$ by what has
already been proved.  Hence $\cl$ meets $S$ in at least two points,
namely $P_1$ and $P_2$, and our hypothesis implies that every point of
$\cl$ belongs to $S$.  In particular $P\in S$, as required.

The proof of part (iii) of Lemma \ref{l2} is the most involved.  
We will write
$S^{c}$ for the complement $\Aft\setminus S$ of $S$.  Our strategy
will be to show that if $S^c$ also contains $t+1$ points in general 
position, then both $\# S>\tfrac12 q^t$ and $\# S^c>\tfrac12 q^t$,
which will provide a contradiction.  We observe that the hypothesis of
part (iii) is symmetric between $S$ and $S^c$, since $S$ meets $\cl$
in at least two points if and only if $\#(S^c\cap\cl)<q-1$.  Now let
$R$ be either $S$ or $S^c$, and assume that $R$ contains $t+1$ points
$P_0,\ldots,P_{t}$ 
in general position.  For our inductive assumption we suppose that,
for any linear space $L\subset\Aft$ of dimension $t-1$, either $R\cap
L$ fails to contain $t$ points in general position, or $L\setminus R$
is contained in a proper linear subspace of $L$.  Thus either $L\cap
R$ or $L\cap R^c$ lies in a proper linear subspace of $L$.  When
$L=L_0$ is generated by $P_1,\ldots,P_t$ we must be in the second case.
For every $P\in L_0$, the line $\cl$ generated by $P$ and $P_0$ meets
$R$ in at least two points (namely $P$ and $P_0$) and hence contains
at least $q-1$ points of $R$.  For distinct choices of $P$ the sets
$\cl\setminus\{P_0\}$ are disjoint, whence
\begin{equation}\label{e1}
\# R\ge 1+(q-2)\#(L_0\cap R).
\end{equation}
Now suppose that every linear space $L$ of dimension $t-1$,
parallel but not equal to $L_0$, has the property that $L\cap R$ lies in
a proper linear subspace of $L$.  Then since $\Aft$ is a disjoint
union of $L_0$ with the various spaces $L$, we see that
\[\# R\le \#(L_0\cap R)+(q-1)q^{t-2}.\]
Comparing this with (\ref{e1}) yields $(q-3)\#(L_0\cap
R)<(q-1)q^{t-2}$, and since $\#(L_0\cap R)\ge q^{t-1}-q^{t-2}$ we find
that $q<4$.  Thus, under the assumption that $q\ge 4$, there must be
at least one space $L_1$, parallel but not equal to $L_0$, for 
which it is $L_1\cap R^c$ which is
contained in a proper linear subspace of $L_1$.  If we pick any point
$Q\in L_1$ and count points of $R$ on lines from $Q$ to $L_0$ we will
obtain at least $(q-2)\#(L_0\cap R)$ points of $R$ not lying on $L_1$,
by the argument that established (\ref{e1}).  Allowing for points of
$L_1\cap R$ we find that
\[\# R\ge (q-2)\#(L_0\cap R)+\#(L_1\cap R).\]
However we have arranged that $L_0\cap R^c$ and $L_1\cap R^c$ are both
contained in proper linear subspaces, so that $\#(L_0\cap R)\ge
q^{t-1}-q^{t-2}$, and similarly for $\#(L_1\cap R)$.  It then follows
that $\# R\ge (q-1)^2q^{t-2}>\tfrac12 q^t$, since $q\ge 4$.  As
explained above, this inequality leads to the claim made in part (iii)
of the lemma.

Finally we turn to part (iv) of Lemma \ref{l2}.
For the induction step we assume that we have 
a linear space $L_0$ of dimension $t-1$, along with a point $P_0$ not in $L_0$,
such that 
\[\#(S\cap L_0)\ge\frac{m^{t}-1}{m-1}\]
and $P_0\in S$.  If $P$ is any point in $S\cap L_0$ the line $\cl$
generated by $P$ and $P_0$ contains at least 2 points of $S$, and
hence by our hypothesis contains at least $m+1$ such points.  For
different points $P$ the sets $\cl\setminus \{P_0\}$ are disjoint,
whence
\[\# S\ge 1+m\#(S\cap L_0)\ge 1+m\frac{m^{t}-1}{m-1}=\frac{m^{t+1}-1}{m-1}.\]
This completes the induction.
\bigskip

We now establish Theorem \ref{geom}.  If any of the polynomials $\Bf$
vanishes identically we may clearly remove it from the collection, and
if any of the polynomials is a non-zero constant then $z(\Bf;\Af)$ is
empty, so that the result is trivial.  It follows that we may assume
that each of the $f_i$ has positive degree, whence $r$, which is the 
number of polynomials $f_i$, is at most $d$.  Thus the zero-set of the
polynomials $\Bf$, over $\overline{\mathbb{F}_q}$ is an algebraic
variety, $V$ say, of dimension at most $n$.  This
variety is defined over $\F_q$, but is not necessarily absolutely 
irreducible.  We shall call an irreducible component of $V$ ``bad'' if
it is not defined over $\F_q$.  Thus our
goal is to replace $V$ by a variety $V^*$ with no bad components, and 
such that $V(\F_q)=V^*(\F_q)$.  We proceed to
describe a general ``reduction process'' for our varieties.  In order
to keep track of the number of absolutely irreducible components that
these varieties have, and of their degrees, it will be convenient to
define $\delta(X)=\sum\deg(X_i)$ whenever the variety $X$ is a union
of absolutely irreducible components $X_i$.  It follows that
$\delta(V)\le C_0(n,d)$ for some number $C_0(n,d)$ depending on $n$
and $d$ alone, and we aim to establish an analogous bound
$\delta(V^*)\le C(n,d)$ for the variety $V^*$.

Suppose
that $V$ has at least one bad component. Choose such a component, 
$W$ say, of maximal
dimension $D$ say. Let the Galois conjugates of $W$ be $W=W_1,\ldots,W_k$.
Since each of these is an irreducible component of $V$ we see that
$k\le C_0(n,d)$.  Any $\F_q$-point of $W$ lies on
$W':=W_1\cap\ldots\cap W_k$.  Thus if we replace the components
$W_1,\ldots,W_k$ of $V$ by $W'$, and call the resulting variety $V_1$,
we see that $V(\F_q)=V_1(\F_q)$.  Moreover $V_1$ contains no bad
components of dimension greater than $D$, and since $\dim(W')<\dim(W)$
there is one fewer bad component of dimension $D$.  We also note that
$\delta(W')$ can be bounded in terms of $n$ and $C_0(n,d)$ alone
whence $\delta(V_1)\le C_1(n,d)$ for a suitable integer $C_1(n,d)$.
We now repeat the reduction process, passing from $V_1$ to $V_2$ and so
on.  Since the number of bad components of maximal dimension is
reduced at each step, while $\delta(V_i)$ remains under control, the
process will eventually terminate, and will produce a variety $V^*$
for which $\delta(V^*)$ is bounded in terms of $n$ and $d$.  The theorem
then follows.

\bigskip
\bigskip

Mathematical Institute,

Radcliffe Observatory Quarter

Woodstock Road

Oxford

OX2 6GG

UK
\bigskip

{\tt rhb@maths.ox.ac.uk}

\end{document}